\documentclass[10pt]{amsart}

\newif\ifamsfonts
\amsfontstrue
\ifamsfonts
\usepackage{amsmath,amsthm,amssymb}
\newcommand{\N}{{\mathbb N}}
\newcommand{\Z}{{\mathbb Z}}

\else
\newcommand{\N}{{\bf N}}
\newcommand{\Z}{{\bf Z}}

\fi

\renewcommand{\text}[1]{\mbox{#1}}

\newcommand{\inn}{\mbox{int}\,}

\newtheorem {lemma}{Lemma}
\newtheorem {theo}{Theorem}

\newtheorem {prop}{Proposition}

\newcommand{\eps}{\varepsilon}
\newcommand{\labda}{\lambda}

\newcommand{\pr}{\noindent{\em Proof: }}
\newcommand{\diam}{\mbox{diam}}
\newcommand{\orb}{\mbox{orb}}
\renewcommand{\a}{\mbox{P}}
\newcommand{\T}{{\mathcal T}}
\newcommand{\ie}{{\em i.e. }}
\newcommand{\eg}{{\em e.g. }}
\newcommand{\para}{\newline \indent}

\begin{document}
\title{Homeomorphic Restrictions of Unimodal Maps}
\author{Henk Bruin} 
\address{
Department of Mathematics,
California Institute of Technology,
Pasadena CA 91125 USA}
\email{bruin@cco.caltech.edu}

\date{April 5 1999.}

\begin{abstract}
\noindent
Examples are given of tent maps $\T$ for which there
exist non-trivial sets $B \subset [0,1]$ such that $\T:B \to B$
is a homeomorphism.
\end{abstract}

\maketitle

\section{Introduction}
Let $\T:[0,1] \to [0,1]$ be a unimodal map, \ie
$\T$ is a continuous map with a unique turning point $c \in [0,1]$
such that $\T|_{[0,c]}$ is increasing and $\T|_{[c,1]}$ is decreasing.
Obviously, $\T:[0,1] \to [0,1]$ is not homeomorphic, but we can ask
ourselves if there are sets $B \subset [0,1]$ such that
$\T:B \to B$ is homeomorphic.
If $B$ is a union of periodic
orbits, then this is obviously the case. Also if there is a subinterval
$J \subset [0,1]$ such that $|\T^n(J)| \to 0$, it is easy to construct
an uncountable set $B \subset \cup_{n \in \Z} \T^n(J)$ such that
$\T:B \to B$ is homeomorphic.
A third example is $B = \omega(c)$
($\omega(x) := \cap_i \overline{\cup_{j \geq i} \T^j(x)}$),
when $\T$ is infinitely renormalizable.
Indeed, in this case $\omega(c)$ is a so-called solenoidal attractor and
$\T:\omega(c) \to \omega(c)$ is topological conjugate to an adding machine
and therefore a homeomorphism.
For these and other general results on unimodal maps, see
\eg \cite{ALM,dMvS}.
\para
In \cite{BOT} the above question was first raised, and properties
of $B$ were discussed. To avoid the mentioned trivial examples
let us restrict the question for maps $\T$ that are {\em locally eventually
onto}, \ie every interval $J \subset [0,1]$ iterates
to large scale: $\T^n(J) \supset [\T^2(c), \T(c)]$ for $n$ sufficiently
large. Because every locally eventually onto unimodal map
is topologically conjugate to a some tent map $T_a$, $T_a(x) =
\min(ax,a(1-x))$ with $a > \sqrt2$, we can restate the question to
\begin{quote}
Are there tent maps $T_a$, $a > \sqrt2$, 
that admit an infinite compact set $B$
such that $T_a:B \to B$ is a homeomorphism?
\end{quote}
\noindent
This turns out to be the case.
To be precise, we prove
\begin{theo}\label{main}
There exists a locally uncountable dense set $A \subset [\sqrt{2},2]$
such that $T_a:\omega(c) \to \omega(c)$ is homeomorphic
for every $a \in A$.
\end{theo}
\noindent
We remark that because $\omega(c)$ is nowhere dense for each $a \in A$,
$A$ is a first category set of zero Lebesgue measure,
see \cite{BDOT,BM}.

This paper is organized as follows.
In the next section, we discuss, following \cite{BOT}, some properties
that $B$ has to satisfy. In section 3 we recall some facts from kneading
theory. Theorem~\ref{main} is proven in section 4 and in the last section 
we give a different construction to solve the main question.

{\bf Acknowledgement:} I would like to thank the referee for the careful
reading and valuable remarks.

\section{Properties of $B$}

Throughout the paper we assume that $\T = T_a$ is a tent map with slope
$a > 1$ and that $B$ is a compact infinite set such that
$\T:B \to B$ is a homeomorphism.

\begin{prop}\label{propB}
Under the above assumptions, $B = \omega(c)$ modulo a countable set, and
$\omega(c)$ is minimal.
\end{prop}
\noindent
Let us first recall a result of Gottschalk and Hedlund \cite{GH}.
A self map $f$ on a compact metric space is called {\em locally expanding}
if there exist $\eps_0 > 0$ and $\labda > 1$ such that
$d(f(x),f(y)) > \labda d(x,y)$ whenever $d(x,y) < \eps_0$.

\begin{lemma}[\cite{GH}]\label{gothed}
If $X$ is a compact metric space and $f:X \to X$ is a locally expanding
homeomorphism, then $X$ is finite.
\end{lemma}

\pr
Because $f^{-1}$ is continuous and $X$ is compact, there exists $\delta > 0$
such that $d(x,y) < \eps$ implies $d(f^{-1}(x), f^{-1}(y)) < \delta$.
Obviously $\delta$ can be taken small as $\eps \to 0$.
In particular, if $\eps \ll \eps_0$,
local expandingness gives that we can take $\delta = \frac1{\labda} \eps$.
Let $\cup_i U_i$ be an open cover of $X$ such that $\diam(U_i) < \eps$ for
each $i$. As $X$ is compact we can take a finite
subcover $\cup_{i=1}^N U_i$.
By definition of $\delta$, $\diam(f^{-1}U_i) < \delta$ and as $f$ is locally
expanding, $\diam(f^{-1}U_i) < \frac1{\labda} \diam(U_i) <
\frac{\eps}{\labda}$.
Repeating this argument, we obtain for each $n$ a finite cover
$\cup_{i=1}^N f^{-n}U_i$ of $X$ and
$\diam(f^{-n}U_i) < \labda^{-n}\eps \to 0$ uniformly.
Hence $X$ must be finite.
\qed

\medskip
\noindent
The proof of Proposition~\ref{propB} uses ideas from \cite{BOT}:
\medskip

\pr
The map $\T$ is locally expanding on every compact set that excludes $c$.
Therefore, if $c \notin B$, the previous lemma shows that $B$ is finite.
Assume therefore that $c \in B$, and hence $\omega(c) \subset B$.
If $\omega(c)$ is not minimal, then there exists $x \in \omega(c)$
such that $c \notin \omega(x)$. Then $\T: \omega(x) \to \omega(x)$ is a
homeomorphism. Hence by Lemma~\ref{gothed}, $x$ must be a periodic
orbit, say with period $N$.
Take $U \owns x$ so small that for each $0 \leq i < N$,
$\T^{-1}(\T^i(U))$ has only one component that intersects $B$.
This is possible because $\T:B \to B$ is one-to-one and $c \notin \orb(x)$.
Let $m$ be minimal such that $c_m \in \cup_{i=0}^{N-1} \T^i(U)$,
say $c_m \in \T^i(U)$.
But then $c_{m-1} \in B$ belongs to a component of $\T^{-1}(\T^i(U))$
that does not intersect $B$. This contradiction shows that $\omega(c)$ is
minimal.
\para
Now assume that $x \in B \setminus \omega(c)$ is such that $c \in
\omega(x)$.
Let $0 < \eps < d(x,\omega(c))$, and let $U_1 = B(c;\frac{\eps}{2})$ be
the open $\frac{\eps}{2}$-ball centered at $c$.
By taking $\eps$ smaller if necessary we can assume that if $U$ is any
interval disjoint from $\T(U_1)$ and with $\diam(U) < \eps$,
at most one component of $\T^{-1}(U)$ intersects $B$.
Finally assume that 
$\T^n(\partial U_1) \cap U_1 = \emptyset$ for all $n \geq 1$.
This happens if \eg $\partial U_1$ contains the point in a
periodic orbit which is closest to the critical point. Because
$c$ is an accumulation point of such periodic points 
(for tent maps $T_a$ with $a > 1$), this last assumption can be realized.
\para
For $i \geq 1$ define $U_{i+1}$ to be the component of $\T^{-1}(U_i)$ that
intersects $B$. As $\diam(U_{i-1}) < \eps$ and $\T$ has slope $> 1$,
$\diam(U_{i+1}) < \eps$ and we can continue the construction,
at least as long as $U_i \cap \T(U_1) = \emptyset$.
Let $N$ be minimal such that $U_N \cap U_1 \neq \emptyset$.
Because $c$ is recurrent, $N$ exists.
Then $U_N \subset U_1$, because otherwise
$\partial U_1 \subset U_N$ and $\T^N(\partial U_1) \subset U_1$.
This would contradict the assumption on $\partial U_1$.
\para
Because $c \in \omega(x)$, there exists $m$ minimal such that
$\T^m(x) \in \cup_{i=1}^N U_i$, say $\T^m(x) \in U_i$.
Then, as before, $\T^{m-1}(x) \in B$ lies in a component of $\T^{-1}(U_1)$ that
is disjoint from $B$. This contradiction shows $c \notin \omega(x)$ and
using the above arguments $x$ must be a periodic point. Therefore
$B \subset \omega(c)$ up to a countable set.
\qed

\section{Preliminaries about Kneading Theory}

Let us start with some combinatorics of unimodal maps.
Write $c_n := \T^n(c)$.
We define {\em cutting times} and the {\em kneading map}
of a unimodal map. These ideas were introduced by Hofbauer, see
\eg \cite{H}. A survey can be found in \cite{bruin}.
\para
If $J$ is a maximal (closed) interval on which $\T^n$ is monotone, then
$\T^n:J \to \T^n(J)$ is called a {\em branch}. If $c \in \partial J$,
$\T^n:J \to \T^n(J)$ is a {\em central branch}.
Obviously $\T^n$ has two central branches, and they have the same
image. Denote this image by $D_n$.
\para
If $D_n \owns c$, then $n$ is called a {\em cutting time}.
Denote the cutting times by $\{ S_i\}_{i \geq 0}$,
$S_0 < S_1 < S_2 < \dots$ For interesting
unimodal maps (tent maps with slope $> 1$) $S_0 = 1$ and $S_1 = 2$.
The sequence of cutting times completely determines the tent map and vice
versa.
It can be shown that $S_k \leq 2S_{k-1}$ for all $k$. Furthermore, the
difference between two consecutive cutting times is again a cutting time.
Therefore we can write
\begin{equation}\label{recur_cut}
S_k = S_{k-1} + S_{Q(k)},
\end{equation}
for some integer function $Q$, called the {\em kneading map}.
Each unimodal map therefore is characterized by its kneading map.
Conversely, each map $Q:\N \to \N \cup \{ 0 \}$ satisfying
$Q(k) < k$ and the {\em admissibility condition}
\begin{equation}\label{admis}
\{ Q(k+j)\}_{j \geq 1} \succeq \{ Q(Q^2(k)+j) \}_{j \geq 1}
\end{equation}
(where $\succeq$ denotes the lexicographical ordering)
is the kneading map of some unimodal map.
Using cutting times and kneading map, the following properties
of the intervals $D_n$ are easy to derive:
$$
D_{n+1} = \left\{
\begin{array}{ll}
\T(D_n) & \text{ if } c \notin D_n, \\ \relax
[c_{n+1},c_1] & \text{ if } c \in D_n.
\end{array} \right.
$$
Equivalently:
\begin{equation}\label{Dn}
D_n = [c_n,c_{n-S_k}], \text{ where } k = \max\{ i; S_i < n\},
\end{equation}
and in particular
$$D_{S_k} = [c_{S_k}, c_{S_{Q(k)}}].$$
Let $z_k < c < \hat z_k$ be the boundary points of the domains the two central branches of
$\T^{S_{k+1}}$.
Then $z_k$ and $\hat z_k$ lie in the interiors of the domains of the central branches
of $\T^{S_k}$ and  $\T^{S_k}(z_k) = \T^{S_k}(\hat z_k) = c$. Furthermore,
$\T^j$ is monotone on $(z_k,c)$ and $(c, \hat z_k)$ for all $0 \leq j \leq S_k$.
These points are called {\em closest precritical points}, and the
relation (\ref{recur_cut}) implies
\begin{equation}\label{pos_cSk}
\T^{S_{k-1}}(c) \in (z_{Q(k)-1}, z_{Q(k)}] \cup [\hat z_{Q(k)}, \hat z_{Q(k)-1}).
\end{equation}
We will use these relations repeatedly without specific reference.

Let us also mention some relations with the standard kneading theory.
The {\em kneading invariant} $\kappa = \{ \kappa_n \}_{n \geq 1}$
is defined as
$$
\kappa_n = \left\{
\begin{array}{ll}
0 & \text{ if } \T^n(c) < c, \\
\ast & \text{ if } \T^n(c) = c, \\
1 & \text{ if } \T^n(c) > c.
\end{array} \right.
$$
If we define
\begin{equation}\label{tau}
\tau: \N \to \N, \,\,\, \tau(n) = \min\{ m > 0; \kappa_m \neq
\kappa_{m+n} \},
\end{equation}
then we retrieve the cutting times as follows:
$$
S_0 = 1 \text{ and } S_{k+1} = S_k + \tau(S_k)
= S_k + S_{Q(k+1)} \text{ for } k \geq 0.
$$
In other words (writing $\kappa_i' = 0$ if $\kappa_i = 1$ and vice versa),
\begin{equation}\label{knea}
\kappa_1 \dots \kappa_{S_k} =
\kappa_1 \dots \kappa_{S_{k-1}}\kappa_1 \dots \kappa'_{S_{Q(k)}}.
\end{equation}
For the proofs of these statements, we refer to \cite{bruin}.

\section{Proof of the Theorem}

Let us introduce an adding machine-like
number system that factorizes over the action $\T:\omega(c) \to \omega(c)$:
Let $\{ S_k \}$ be the cutting times of a unimodal map
and assume that the corresponding kneading map $Q$ tends to infinity.
Any non-negative integer $n$ can be written in a canonical way as
a sum of cutting times: $n = \sum_i e_i S_i$,
where 
$$e_i = \left\{
\begin{array}{ll}
1 &\text{ if } i = \max\{ j; S_j \leq n- \sum_{k > i} e_k S_k \},\\
0 &\text{ otherwise.}
\end{array} \right.
$$
In particular $e_i = 0$ if $S_i > n$.
In this way we can code the non-negative integers $\N$ as zero-one sequences
with a finite number of ones: $n \mapsto \langle n \rangle \in \{ 0,1 \}^{\N}$.
Let $E_0 = \langle \N \rangle$ be the set such sequence, 
and let $E$ be the closure of $E_0$ in the product topology.
This results in
$$
E= \{ e \in \{ 0,1 \}^{\N}; e_i = 1 \Rightarrow e_j = 0 \text{ for }
Q(i+1) \leq j < i \},
$$
because if $e_i = e_{Q(i+1)} = 1$, then this should be
rewritten to $e_i  = e_{Q(i+1)} = 0$ and $e_{i+1} = 1$.
 It follows immediately that
for each $e \in E$ and $j \geq 0$,
\begin{equation}\label{rule}
e_0S_0 + e_1 S_1 + \dots + e_j S_j < S_{j+1}.
\end{equation}

There exists the standard addition of $1$ by means of `add and carry'.
Denote this action by $\a$. Obviously
$\a(\langle n \rangle) = \langle n+1 \rangle$. 
It is known (see \eg \cite{BKS,GLT}) that $\a:E \to E$ is continuous if and only if $Q(k) \to \infty$,
and that $\a$ is invertible on $E \setminus \{ \langle 0 \rangle \}$.
The next lemma describes the inverses of $\langle 0 \rangle$ precisely.

\begin{lemma}\label{-1}
For a sequence $e \in E$, let $\{ q_j\}_{j \geq 0}$ be the index set (in increasing order)
such that $e_{q_j} = 1$. We have $\a(e) = \langle 0 \rangle$ if and only if
$e \notin E_0$, $Q(q_0+1) = 0$ and $Q(q_j+1) = q_{j-1}+1$ for $j \geq 1$.
\end{lemma}

\pr
This follows immediately from the add and carry construction, because the
condition on $\{ q_j\}$ is the only way the addition of $1$ carries 
`to infinity'.
\qed

\medskip

\noindent
The next lemma gives conditions under which $\a$ is invertible
on the whole of $E$.

\begin{lemma}\label{invertible} Let $Q$ be a kneading map such that $Q(k) \to \infty$.
Suppose that there is an infinite sequence $\{ k_i \}$
such that for all $i$ and $k > k_i$ 
\begin{itemize}
\item either $Q(k) \geq k_i$, 
\item or $Q(k) < k_i$ and there are only finitely many $l > k$ such that
$Q^n(l) = k$ for some $n \in \N$, 
\end{itemize}
then $\a$ is a homeomorphism of $E$.
\end{lemma}

\pr
Because $\a$ is continuous and invertible outside $\langle 0 \rangle$
and $E$ is compact, it suffices to show that 
$\# \a^{-1}(\langle 0 \rangle) = 1$.
Let $e \in \a^{-1}(\langle 0 \rangle)$ and let $\{ q_j \}_j$ be the index sequence of
the non-zero entries of $e$.
By the previous lemma and our assumption, we see that
$\{ k_i - 1 \}$ must be a subsequence of $\{ q_j\}_j$.
But because any $q_j$ determines $q_{j'}$ for $j' < j$, there can only
be one such sequence $\{ q_j\}_j$ and one preimage $e$.
\qed

\medskip
\noindent
{\bf Example 1} Take any sequence $\{ k_i\}$ with $k_i > k_{i-1}+10$
and define $Q$ as
$$
Q(k) = \left\{ \begin{array}{ll}
k-4 & \text{ if } k > K \text{ and } k-4 \in \{ k_i\}, \\
k-3 & \text{ if } k > K \text{ and } k-5 \in \{ k_i\}, \\
\text{arbitrary} & \text{ if } k \leq K, \\
k-2 &\text{ otherwise,}
\end{array} \right.
$$
provided $Q$ is admissible. 
It is shown in \cite{bruin} that $Q$ belongs to a renormalizable
map of period $S_k$ if and only if 
\begin{equation}\label{renor}
Q(k+1) = k \text{ and } Q(k+j) \geq k \text{ for all } j \geq 1. 
\end{equation}
Because $Q(k) \leq k-2$ for $k \geq K$, we can avoid renormalizable maps.
This shows that there is a locally uncountable
dense set of tent maps, whose kneading maps $Q$ satisfy
Lemma~\ref{invertible} and $Q(k) \to \infty$. In fact, this example 
also satisfies 
conditions (\ref{strong_admis}) and (\ref{stop_carry}) 
of Theorem~\ref{homeo} below.

\medskip

\noindent
{\em Remark:}  Lemmas~\ref{-1} and \ref{invertible} 
also indicate how to construct 
a number systems $(E,P)$ where for $d \in \N \cup \{ \aleph_0\}$,
$\langle 0 \rangle$ has exactly $d$ preimages. 
For example, if $Q(k) = \max(0, k-d)$, then $\langle 0 \rangle$
is $d$ preimages.
The rest of this section makes 
clear that this yields examples of maps where $\T : \omega(c) \to \omega(c)$ 
is one-to-one except for $d$ points in $\cup_{n \geq 0} \T^{-n}(c)$.

\medskip

Given $n \in (S_{k-1}, S_k]$, define $\beta(n) = n-S_{k-1}$. 
It is easy to check that $\langle \beta(n) \rangle$ is
$\langle n \rangle$ with the last non-zero entry changed to $0$.
The map $\beta$ also has a geometric interpretation in the Hofbauer tower:
It was shown in \cite[Lemma 5]{BKS} that for all $n \geq 2$,
\begin{equation} \label{nest1}
D_n \subset D_{\beta(n)}.
\end{equation}
In fact, $D_n$ and $D_{\beta(n)}$ have the boundary 
point $c_{\beta(n)}$ in common.
Recall that for $e \in E$, $\{ q_j\}_j$ is the index sequence 
of the non-zero entries of $e$.
Define
$$
b(i) := \sum_{j \leq q_i} e_j S_j.
$$
We have $b(i) \geq S_{q_i}$ by definition of $q_i$ and
$b(i) < S_{q_i+1}$ by (\ref{rule}).
It follows that $\beta(b(i)) = b(i-1)$.
By a {\em nest} of levels will be meant a sequence of levels $D_{b(i)}$.
By (\ref{nest1}) and the fact that $\beta(b(i)) = b(i-1)$, 
these levels lie indeed nested, and because  
$Q(k) \to \infty$ implies that $|D_n| \to 0$ (see \cite{bruin}), 
each nest defines
a unique point $x = \cap_i D_{b(i)} \in \omega(c)$.
Therefore the following projection (see \cite{BKS}) makes sense:
$$\pi(\langle n \rangle) = c_n$$
and
\begin{equation}\label{nest}
\pi(e \notin E_0) = \cap_i D_{b(i)}.
\end{equation}
Obviously $\T \circ \pi = \pi \circ \a$ and it can be shown
(see \cite[Theorem 1]{BKS}) that $\pi:E \to E$ is uniformly continuous
and onto. 

\noindent
Note that a nest contains exactly one cutting level $D_{b(0)}$.
If $\{ D_{b(i)} \}_i$ is some nest converging to $x$, then $\{ \T(D_{b(i)}) \}_i$ 
is a nested sequence of levels converging to $\T(x)$. To obtain a nest 
of $\T(x)$, we may have to add or delete
some levels, but $\{ \T(D_{b(i)}) \}$ asymptotically coincides with a nest converging to
$\T(x)$.

\begin{theo}\label{homeo}
If $Q$ is a kneading map that satisfies Lemma~\ref{invertible} as well as
\begin{equation}\label{strong_admis}
Q(k+1) > Q(Q^2(k)+1) + 1 
\end{equation}
for all  $k$ sufficiently large, and 
\begin{equation}\label{stop_carry}
Q(s+1) = Q(\tilde s+1) \text{ for } s \neq \tilde s \text{ implies } 
Q^{n+1}(s) \neq Q^{\tilde n+1}(\tilde s),
\end{equation}
for any $n, \tilde n \geq 0$ such that 
$Q^n(s) \neq Q^{\tilde n}(\tilde s)$,
then any map $\T$ with kneading map $Q$ is homeomorphic on $\omega(c)$.
\end{theo}

\pr
In view of Lemma~\ref{invertible} we only 
have to show that $\pi:E \to \omega(c)$ is one-to-one.
First note (see also \cite[Theorem 1]{BKS}) that $\pi^{-1}(c) = \langle 0 \rangle$.
Indeed, if $e \neq \langle 0 \rangle$ and $\pi(e) = c$, then, taking $k < l$ the
first non-zero entries of $e$, $c \in D_{S_k+S_l}$. Then $S_k + S_l$ is 
a cutting time $S_m$ and we have $m = l+1$ and $k = Q(l+1)$. 
This would trigger a carry to $e_k = e_l = 0$ and $e_m = 1$.
Because $P$ is invertible, also $\# \T^{-n}(c) \cap \omega(c) = 1$
for each $n \geq 0$.
Assume from now on that $x \in \omega(c) \setminus \cup_{n \geq 0} \T^{-n}(c)$.
We need one more lemma:

\begin{lemma}\label{overlap}
Let $\T$ be a unimodal map whose kneading map satisfies (\ref{strong_admis})
and tends to infinity.
Then there exists $K$ such that for any $n \notin \{ S_i \}_i$ such that 
$\beta(n)$ is a cutting time (\ie $n = S_r+S_t$ for some $r < t$
with $r < Q(t+1)$),
and every $k \geq K$,
$\inn D_n$ does not contain both $c_{S_k}$ and a point 
from $\{ z_{Q(k+1)-1}, \hat z_{Q(k+1)-1}\}$.
\end{lemma}

\pr
Assume the contrary. 
Write $n = S_r + S_t$ with $r < Q(t+1)$ and let $k$ but such that
$z_{Q(k+1)-1}$ or $\hat z_{Q(k+1)-1} \in D_n \subset D_{S_r}$.
Formula (\ref{pos_cSk}) implies that $Q(r+1) < Q(k+1)$, see figure 1.

\begin{center}
\unitlength=8mm
\begin{picture}(12,4.3)
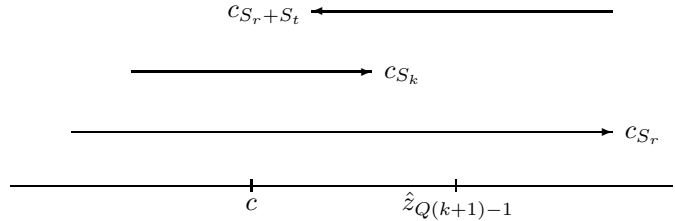
\put(2,-0.3){Figure 1: The levels $D_{S_k}$ and $D_{S_r+S_t}$}

\put(0, 1.1){\line(1,0){11}} 
\put(4, 1){\line(0,1){0.2}} \put(3.9, 0.7){$c$} 
\put(7.4, 1){\line(0,1){0.2}} \put(6.5, 0.7){$\hat z_{Q(k+1)-1}$} 

\put(1,2){\vector(1,0){9}} \put(10.2, 1.9){$c_{S_r}$} 
\put(2,3){\vector(1,0){4}} \put(6.2, 2.9){$c_{S_k}$} 
\put(10,4){\vector(-1,0){5}} \put(3.6, 3.9){$c_{S_r+S_t}$} 
\end{picture}
\end{center}

\bigskip
\noindent
It follows that also $z_{Q(r+1)}$ or $\hat z_{Q(r+1)} \in D_{S_r +S_t}$
and therefore $S_r+S_t+S_{Q(r+1)} = S_{t+1}$.
This gives
\begin{equation}\label{x1}
r+1 = Q(t+1),
\end{equation}
and $S_r + S_t = S_{t+1} - S_{Q^2(t+1)}$.
If also $z_{Q(r+1)+1}$ or $\hat z_{Q(r+1)+1} \in D_{S_r+S_t}$, 
then $S_{Q(r+1)+1} - S_{Q(r+1)} +S_{t+1} = S_{t+2}$,
which yields $Q(Q(r+1)+1) = Q(t+2)$. Using (\ref{x1}) for $t+1$, this gives
$Q(t+1+1) = Q(Q^2(t+1)+1)$.
This contradicts (\ref{strong_admis}),
if $t$ is sufficiently large. For smaller $t$, there are only finitely
many pairs $r < t$. For $k$ sufficiently large
(recall that $Q(k+1) \to \infty$, so $c_{S_k} \to c$), 
$D_{S_r+S_t} \not\owns c_{S_k}$.
Hence
\begin{equation}\label{x2}
D_{S_r+S_t} \text{ contains at most one closest precritical point.}
\end{equation}
Therefore, as $c_{S_k} \in D_{S_r+S_t}$, $Q(r+1) = Q(k+1)-1$.
Take the $S_{Q(r+1)}$-th iterate of $D_{S_r +S_t}$ and $[c,c_{S_k}]$ to obtain
$D_{S_k+S_{Q(r+1)}} \cap D_{S_{t+1}} \neq \emptyset$, see figure 2.

\begin{center}
\unitlength=8mm
\begin{picture}(12,3.3)
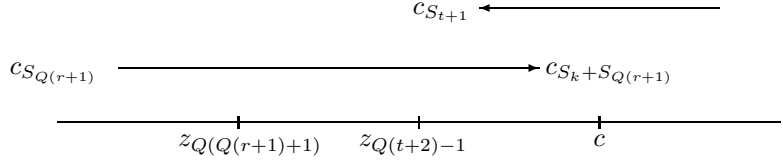
\put(2,-0.3){Figure 2: The levels $D_{S_{Q(r+1)}+S_k}$ and $D_{S_{t+1}}$}

\put(0, 1.1){\line(1,0){12}} 
\put(9, 1){\line(0,1){0.2}} \put(8.9, 0.7){$c$} 
\put(3, 1){\line(0,1){0.2}} \put(2, 0.7){$z_{Q(Q(r+1)+1)}$} 
\put(6, 1){\line(0,1){0.2}} \put(5, 0.7){$z_{Q(t+2)-1}$} 

\put(1,2){\vector(1,0){7}} \put(-0.8, 1.9){$c_{S_{Q(r+1)}}$}  
  \put(8.1, 1.9){$c_{S_k+S_{Q(r+1)}}$}
\put(11,3){\vector(-1,0){4}} \put(5.9, 2.9){$c_{S_{t+1}}$} 
\end{picture}
\end{center}

\bigskip
\noindent
By (\ref{x1}) we have $Q(Q(r+1)+1) = Q(Q^2(t+1)+1)$, and using 
(\ref{strong_admis}) on $t+1$, we obtain
$$
Q(Q^2(t+1)+1) < Q(t+1+1)-1 = Q(t+2)-1.
$$
Hence there are at least two closest precritical points contained in 
$D_{S_k + S_{Q(r+1)}}$. This contradicts the arguments leading to
(\ref{x2}).
\qed

\medskip
\noindent
We continue the proof of Theorem~\ref{homeo}.
Observe that if $\pi(e) = \pi(\tilde e) = x$ for some $e \neq \tilde e$,
then the corresponding nests $\{ D_{b(i)} \}$ and $\{ D_{\tilde b(i)} \}$
are different, but both nests converge to $x$. 
Because $x \notin \cup_{n \geq 0} \T^{-n}(c)$, 
$\# \pi^{-1}(\T^n(x)) > 1$ for all $n \geq 0$.
We will derive a contradiction.

\begin{quote}
Claim 1: We can assume that $b(0) \neq \tilde b(0)$.
\end{quote}

Let $i$ be the smallest integer such that
$b(i) \neq \tilde b(i)$, say $b(i) < \tilde b(i)$.
Then also $q_i < \tilde q_i$.
Let $l = S_{\tilde q_i+1} - \tilde b(i)$.
By (\ref{rule}), $l$ is non-negative.
By the choice of $l$, $(P^l(\tilde e))_j = 0$ for all $j \leq \tilde q_i$,
but because $b(i) < \tilde b(i)$ and $l+b(i) < S_{\tilde q_i+1}$,
there is some $j \leq \tilde q_i$ such that $(P^l(e))_j = 1$.

Replace $x$ by $\T^l(x)$, and the corresponding sequences
$e$ and $\tilde e$ by $P^l(e)$ and $P^l(\tilde e)$.
Then for this new point, $q_0 < \tilde q_0$ and $b(0) < \tilde b(0)$.
This proves Claim 1.

\begin{quote}
Claim 2: We can assume that $Q(q_0+1) \neq Q(\tilde q_0+1)$ 
\end{quote}

Assume that $Q(q_0+1) = Q(\tilde q_0+1) =: r$.
We apply $P^{S_r}$ to $e$ and $\tilde e$. 
Write $s = \min\{ j; (P^{S_r}(e))_j = 1\}$ and $\tilde s =
\min\{ j; (P^{S_r}(\tilde e))_j = 1\}$.
To make sure that Claim 1 still
holds, assume by contradiction that $s = \tilde s$.
Then (using (\ref{rule})),
$\sum_{j=0}^{s-1} e_j S_j = \sum_{j=0}^{s-1} \tilde e_j S_j = S_s-S_r$.
But this would imply that $e_j = \tilde e_j$ for all $j < s$, which
is not the case. Hence $s \neq \tilde s$.

From the add and carry procedure it follows that $e_{j-1} = 1$ for 
$j = Q(s)$, $Q^2(s)$, $\dots$, $q_0+1 = Q^n(s)$ for some $n$,
and similarly
$\tilde e_{j-1} = 1$ for 
$j = Q(\tilde s), Q^2(\tilde s), \dots, \tilde q_0+1 = 
Q^{\tilde n}(\tilde s)$ for some $\tilde n$.
As $Q(q_0+1) = Q(\tilde q_0+1)$,  hypothesis
(\ref{stop_carry}) implies $Q(s+1) \neq Q(\tilde s+1)$.
This proves Claim 2. Replace $x$ by $\T^{S_r}(x)$ and the corresponding
sequences $e$ and $\tilde e$ by
$P^{S_r}(e)$ and $P^{S_r}(\tilde e)$.
\para
Note that we can take $q_0 \neq \tilde q_0$  arbitrarily
large, with say $Q(q_0+1) < Q(\tilde q_0+1)$. 
Then we are in the situation of
Lemma~\ref{overlap}, which 
tells us that $D_{b(1)} \cap D_{\tilde b(0)} = \emptyset$.
Therefore the two nests cannot converge to the same point.
This contradiction concludes the proof.
\qed

\medskip

\noindent
{\em Proof of Theorem~\ref{main}:} Combine the previous theorem with Example 1.
\qed

\section{A Different Construction}

In this section we give a different construction which does not involve
the assumption $Q(k) \to \infty$.
Let $k_1$ be arbitrary and $k_2 = k_1+1$.
Put recursively for $i \geq 3$,
$$k_i = 2k_{i-1} - k_{i-2} + 1, \text{ \ie }
k_i-k_{i-1} = k_{i-1}+ k_{i-2} +1.$$
Define the kneading map as
\begin{equation}\label{def1}
Q(k_i) = k_i-1 \text{ for } i \geq 3,
\end{equation}
and choose $Q(k)$ arbitrary for $k \leq k_2$ so that (\ref{admis}) 
and (\ref{strong_admis}) below
are not violated for $k \leq k_2+2$. 
To finish the definition, let 
\begin{equation}\label{def2}
Q(k_i+j) = Q(k_{i-1}+j-1) \text{ for } i \geq 2, 1 \leq j < k_{i+1}-k_i.
\end{equation}
A direct computation shows that (\ref{def1}) and (\ref{def2}) imply
(\ref{strong_admis}) for all $k > k_2$.
Therefore the construction is compatible with the 
admissibility condition (\ref{admis}).
Moreover, (\ref{def1}) and (\ref{def2}) show that condition (\ref{renor})
is not met for $k \geq k_2$. Therefore $Q$ does not 
belong to a renormalizable map of period $\geq S_{k_2}$.

\begin{theo}
If $\T$ is a unimodal map with the kneading map constructed above, then
$\T: \omega(c) \to \omega(c)$ is a homeomorphism.
\end{theo}

\pr
Write $B = \omega(c)$.
Using the levels $D_n$ of the Hofbauer tower, we will
construct covers of $B$ to show that $\T:B \to B$ is a homeomorphism.
Let $\Delta_i = \cup_{n = S_{k_i-1}+1}^{S_{k_i}} D_n$.
We will use the following claims:
\begin{equation}\label{claim1}
c_{S_{k_i-1}} \in [c_{S_k}, 1-c_{S_k}] \text{ for every } k \leq k_i.
\end{equation}
\begin{equation}\label{claim2}
c_n \notin \inn D_{S_{k_i}}
\text{ for } 0 < n \leq S_{k_i}.
\end{equation}
\begin{equation}\label{claim3}
\Delta_i \text{ consists of disjoint intervals.}
\end{equation}
\begin{equation}\label{claim4}
c_n \in \Delta_i \text{ for } S_{k_i-1} \leq n < S_{k_{i+1}-1}
\end{equation}
\begin{equation}\label{claim5}
\Delta_{i+1} \subset \Delta_i
\end{equation}

\pr
Claim (\ref{claim1}): Recall the function $\tau$ from (\ref{tau}). 
Obviously $\tau(n) > \tau(m)$ implies that $c_n \in (c_m, 1-c_m)$.
By construction and equation (\ref{knea}),
$\tau(S_{k_i-1}) = S_{Q(k_i)} = S_{k_i-1} \geq S_{Q(k+1)}$ for
all $k \leq k_i$. Hence
$c_{S_{k_i-1}} \in [c_{S_k}, 1-c_{S_k}]$ for every $k \leq k_i$,

\medskip

Claim (\ref{claim2}):
By construction $D_{S_{k_i}} = [c_{S_{k_i}}, c_{S_{Q(k_i)}}]
= [c_{S_{k_i}}, c_{S_{k_i-1}}]$.
We have 
$$
G_i := \min\{ \tau(S_{k_i}), \tau(S_{k_i-1}) \} 
= S_{k_{i-1}-1}.
$$
Because $\tau(S_k) < G_i$ for
all $k < k_i-1$, $k \neq k_{i-1}-1$, we obtain
$c_{S_k} \notin D_{S_{k_i}}$ for these values of $k$.
If $k = k_i-1$, then $c_{S_k} \in \partial D_{S_{k_i}}$ and not in the
interior.
With respect to $k_{i-1}-1$, note that by (\ref{def2}),
$Q(k_{i-1}-1) = Q(k_i-1)$, 
so $\kappa_{ S_{k_{i-1}-1} } = \kappa_{S_{k_i-1}}$ and  
$c_{S_{k_{i-1}-1}}$ and $c_{S_{k_i-1}}$ lie on the same side of $c$.
Because also $\tau(S_{k_{i-1}-1}) = S_{Q(k_{i-1})} < S_{Q(k_i)} =
\tau(S_{k_i-1})$, $c_{S_{k_{i-1}-1}} \notin D_{S_{k_i}}$.
\para
It remains to consider non-cutting times $n < S_{k_i}$.
Assume by contradiction that $c_n \in \inn D_{S_{k_i}}$,
\ie $D_n$ intersects $D_{S_{k_i}}$ in a non-trivial interval.
Then also $D_{\beta(n)}$ intersects $D_{S_{k_i}}$ where $\beta$ is as in
(\ref{nest1}).
By taking $\beta^j(n)$ instead of $n$ for some $j \geq 0$, we may
assume that $\beta(n)$ is a cutting time.
In particular, $c_n \in \inn D_{S_{k_i}}$ and
$n = S_k + S_t < S_{k_i}$, where $Q(t+1) > k$.
If $k$ is such that $\tau(S_k) < G_i$, then by (\ref{strong_admis})
and Lemma~\ref{overlap}, $D_n \cap D_{S_{k_i}} = \emptyset$.
If $k = k_i-1$,
then $Q(t+1) > k$ implies $t \geq k_i-1$, contradicting that 
$S_k+S_t < S_{k_i}$.
The last possibility is that
$k = k_{i-1}-1$ and $t = k_{i-1}-1$. The above arguments showed that
$c_{S_{k_i-1}}$ and $c_{S_{k_{i-1}-1}}$ lie on the same side of $c$.
Because $\tau(S_{k_{i-1}-1}) < \tau(S_{k_i-1})$, Lemma~\ref{overlap}
applies after all.  
This proves Claim (\ref{claim2}).

\medskip

Claim (\ref{claim3}):
Suppose by contradiction that $D_m \cap D_n \neq \emptyset$
for some $S_{k_i-1} < m < n \leq S_{k_i}$.
Then also
$\T^{S_{k_i}-n}(D_m) \cap \T^{S_{k_i}-n}(D_n) =
D_{m+S_{k_i}-n} \cap D_{S_{k_i}} \neq \emptyset$.
Because $S_{k_i}+m-n$ is not a cutting time, at least one endpoint
of $D_{m+S_{k_i}-n}$ is contained in $D_{S_{k_i}}$.
This contradicts the previous claim.

\medskip

Claim (\ref{claim4}):
Clearly $c_{S_{k_i-1}} \in [c_{S_{k_i}}, c_{S_{k_i-1}}] =
D_{S_{k_i}} \subset \Delta_i$
and for $S_{k_i-1} < n \leq S_{k_i}$, $c_n \in D_n \subset \Delta_i$
by definition.
So let us consider $n = S_{k_i}+1$.
By construction of $Q$ and (\ref{knea}) we obtain
\begin{equation}\label{m}
\begin{array}{rcl}
m &:=& S_{k_{i+1}-1} - S_{k_i} \\
&=& S_{Q(k_i+1)} + S_{Q(k_i+2)} + \dots + S_{Q(k_{i+1}-1)} \\
&=& S_{Q(k_{i-1})} + S_{Q(k_{i-1}+1)} + \dots + S_{Q(k_i-1)} \\
&=& S_{k_i-1} - S_{k_{i-1}-1} = S_{Q(k_i)} - S_{k_{i-1}-1},
\end{array}
\end{equation}
and
$$
\kappa_{S_{k_i}+1} \dots \kappa_{S_{k_{i+1}-1}} =
\kappa_{S_{k_{i-1}-1}+1} \dots \kappa_{S_{Q(k_i)}} =
\kappa_{S_{k_i-1}+S_{k_{i-1}-1}+1} \dots \kappa'_{S_{k_i}}.
$$
Here we `shifted' the word $\kappa_{S_{k_{i-1}-1}+1} \dots \kappa_{S_{Q(k_i)}}$
over $S_{k_i-1}$ entries and used (\ref{knea}) to obtain the second equality. 
Therefore $c_{S_{k_i}+1}$ lies in the same interval of monotonicity
of $\T^{m-1}$ as 
the level 
$D_{S_{k_i-1}+S_{k_{i-1}-1}+1} =
[c_{S_{k_i-1}+S_{k_{i-1}-1}+1}, c_{S_{k_{i-1}-1}+1}]$.
Furthermore
$\T^{m-1}(c_{S_{k_i}+1}) = c_{S_{k_{i+1}-1}}$ and
\begin{equation}\label{mm}
\T^{m-1}(D_{S_{k_i-1}+S_{k_{i-1}-1}+1}) =
\T^{S_{Q(k_i)}-1}(D_{S_{k_i-1}+1}) = D_{S_{k_i}}.
\end{equation}
Claim (\ref{claim1}) gives
$c_{S_{k_{i+1}-1}} \in D_{S_{k_i}}$.
Therefore 
\begin{equation}\label{in_D}
c_{S_{k_i}+1} \in D_{S_{k_i-1}+S_{k_{i-1}-1}+1},
\end{equation} 
and $c_n \in D_{S_{k_i-1}+S_{k_{i-1}-1}+n-S_{k_i}} \subset \Delta_i$
for all $S_{k_i} < n < S_{k_{i+1}-1}$.

\medskip

Claim (\ref{claim5}):
We need to show that
$$
D_{S_{k_{i+1}-1}+j} \subset \Delta_i
\text{ for }
1 \leq j \leq S_{k_{i+1}} - S_{k_{i+1}-1} = S_{k_{i+1}-1}.
$$
Because $\tau(S_{k_i-1}) < \tau(S_{k_{i+1}-1})$,
$c_{S_{k_{i+1}-1}} \in [c_{S_{k_i-1}},1-c_{S_{k_i-1}}]$.
Hence
$D_{S_{k_{i+1}-1}+j} \subset D_{S_{k_i-1}+j}$ for
$0 < j \leq S_{Q(k_i)} = S_{k_i-1}$.

For $j = S_{Q(k_i)}$, $D_{S_{k_i-1}+j} = D_{S_{k_i}}$, and the above line
shows that $D_{S_{k_i}} \supset D_{S_{k_{i+1}-1}+j}$
and these two intervals have the boundary point
$c_{S_{Q(k_i)}}$ in common.
Because also $Q(k_i) = k_i-1$, we get
$D_{S_{k_{i+1}-1}+j} \subset D_{ S_{k_i-1}+(j-S_{k_i-1}) }$ for
$S_{k_i-1} < j \leq S_{k_i-1} + S_{Q(k_i)} = S_{k_i}$.
\para
By formula (\ref{in_D}), one boundary point 
$c_{S_{k_i}+1} \in \partial D_{S_{k_{i+1}-1}+S_{k_i}+1}$ 
belongs to $D_{S_{k_i-1} + S_{k_{i-1}-1}+1}$.
A fortiori,
$D_{S_{k_{i+1}-1}+j} \cap D_{S_{k_i-1}+S_{k_{i-1}-1}+(j-S_{k_i})}
\neq \emptyset$ for
$S_{k_i} < j \leq S_{k_i} + (S_{k_i} - (S_{k_i-1} + S_{k_{i-1}-1}))$.
In particular (cf. (\ref{mm})), for 
$j =  S_{k_i} + (S_{k_i} - (S_{k_i-1} + S_{k_{i-1}-1})) =
S_{k_i} + S_{k_i-1} - S_{k_{i-1}-1}$,
$D_{S_{k_{i+1}-1}+j}$
intersects the level
$D_{ S_{k_i-1} + S_{k_{i-1}-1} + (j - S_{k_i})} = D_{2S_{k_i-1}} =
D_{S_{k_i}}$. (Here we used $Q(k_i) = k_i-1$, \ie 
$S_{k_i} = 2S_{k_i-1}$).
At the same time, by (\ref{m}),
$j = S_{k_i} + S_{k_i-1} - S_{k_{i-1}-1} = S_{k_{i+1}-1}$ and therefore
$D_{S_{k_{i-1}-1} + j} = D_{S_{k_{i+1}}}$.
Thus the intersection is actually an inclusion: 
$D_{S_{k_{i+1}}} \subset D_{k_i}$ and
$D_{S_{k_{i+1}-1}+j} \subset D_{S_{k_i-1}+S_{k_{i-1}-1}+(j-S_{k_i})}$
for all $j$, $S_{k_i} < j \leq S_{k_{i+1}-1}$.
This proves Claim (\ref{claim5}).

\medskip

Let $i$ be arbitrary.
By construction, $\Delta_i \supset \{ c, c_1, \dots, c_{S_{k_i}} \}$.
Claim (\ref{claim5}) used repeatedly gives 
$\orb(c) \subset \Delta_i$,
and because $\Delta_i$ is closed, $B \subset \cap_i \Delta_i$.
Finally, to prove that $\T:B \to B$ is homeomorphic, it suffices to show
that $\T:B \to B$ is one-to-one.
Suppose by contradiction that there exist $y, y' \in B$, $y \neq y'$,
such that $\T(y) = \T(y')$.
Take $i$ so large that
$y$ and $y'$ lie in different intervals of $\Delta_i$.
Say $y \in D_n$ and $y' \in D_m$.
Because $y \neq c \neq y'$, we can assume that
$S_{k_i-1} < m < n < S_{k_i}$.
But then $\T(D_m) \cap \T(D_n) = D_{m+1} \cap D_{n+1} \neq \emptyset$,
contradicting Claim (\ref{claim3}).
This concludes the proof.
\qed

\end{document}